\documentclass{amsart}
\newtheorem{theorem}{Theorem}[section]
\newtheorem{definition}[theorem]{Definition}
\newtheorem{proposition}[theorem]{Proposition}
\newtheorem{lemma}[theorem]{Lemma}
\newtheorem{remark}[theorem]{Remark}

\newtheorem{conjecture}[theorem]{Conjecture}
\newtheorem{corollary}[theorem]{Corollary}

\newcommand{\ZZ}{{\mathbb Z}}

\newcommand{\CC}{{\mathbb C}}

\newcommand{\one}{{\textbf{1}}}
\newcommand{\cP}{{\mathcal P}}
\newcommand{\cQ}{{\mathcal Q}}

\newcommand{\fB}{{\mathfrak B}}

\newcommand{\tr}{\rm{tr} }

\title[ If they are limit periodic? 
]{ If they are limit periodic? 
    }

\author[J. Bellissard, J. Geronimo, A. Volberg, P.
Yuditskii]
{J. Bellissard$^1$, J. Geronimo$^2$, A. Volberg$^3$, P.
Yuditskii$^4$}

\thanks{This work supported by 
$^1$the NSF grant 0300398; 
$^2$the NSF grant DMS -0200219;
$^3$the NSF grant DMS -0200713;
$^4$the Austrian
Founds FWF, project number: P16390--N04}

\address{
School of Mathematics, 
Georgia Institute of Technology, Atlanta, GA
30332-0160}
\email{jeanbel@math.gatech.edu}

\address{
School of Mathematics, 
Georgia Institute of Technology, Atlanta, GA
30332-0160}
\email{geronimo@math.gatech.edu}

\address{Department of Mathematics,
Michigan State University,
East Lansing, MI 48824}
\email{volberg@math.msu.edu}

\address{Institute for Analysis,
Johannes Kepler University of Linz, 
A-4040 Linz, Austria}
\email{Petro.Yudytskiy@jku.at}

\begin{document}

\maketitle

\begin{abstract}
We prove a partial result concerning
the long--standing problem on limit periodicity of the Jacobi matrix 
associated with the balanced measure on the Julia set of an expending
polynomial. Besides this, connections of the problem with
the  Faybusovich--Gekhtman flow and
many other objects
 (the Hilbert transform, the Schwarz derivative,
 the Ruelle and Laplace operators)
 that, we sure,  are of independent
interest, are discussed.
\end{abstract}

\section{Introduction}
In 80's the following interesting 
phenomena was discovered: the spectral measure of an almost periodic
Jacobi matrix can be singular continuous (supported on a Cantor type
set of the zero Lebesgue measure). The effect was studied from both
sides --- from coefficient sequences to spectral data 
\cite {AS}, \cite {BS} and from spectral
data to Jacobi matrices. 

The second, usually more elegant, approach produced 
the following example 
\cite {BBM}, \cite {BGH}.
Let
$T(z)=z^2-C$. For
$C>2$ the Julia set $E$ of $T$ is a real Cantor type
set, $|E|=0$. Denote by
$\mu$ the balanced measure on $E$,
$\mu(T^{-1}(F))=\mu(F)$ for all $F\subset E$. Let
\begin{equation}
J=
\begin{bmatrix}
q_{0}&p_{1}& & & \\
p_{1}&q_{1}&p_{2}& & \\
    &\ddots&\ddots&\ddots & \\
\end{bmatrix}
\end{equation}
be the Jacobi matrix associated to the given measure.
 Note that to construct $J:l^2(\ZZ_+)\to l^2(\ZZ_+)$ one uses the three
term recurrent relation for  polynomials orthonormal in $L^2_{d\mu}$
\begin{equation}\label{r1}
   \lambda P_k(\lambda)=p_k P_{k-1}(\lambda)+q_k
P_{k}(\lambda)+ p_{k+1} P_{k+1}(\lambda),
\end{equation}
of course $P_k\mapsto |k\rangle$, where
$\{|k\rangle\}$ is the
standard basis in $l^2(\ZZ_+)$.

Then the given matrix satisfies the renormalization equation:
$$
V^*T(J)V=J, 
$$ 
where
$V |k\rangle=|2 k\rangle$.
 In fact, this is a system
of nonlinear equations for $p_n$'s ($q_n=0$ in this case), due to which 
at least for $C > 3$ one
gets inductively that 
$$
|p_{ 2^n l+m}-p_m|\le \epsilon_n,
\quad{\rm for\ all} \ \ l, m;\quad \epsilon_n\to 0 \ (n\to \infty).
$$
That is the sequence $\{p_n\}$ and, by  definition
the matrix itself, is limit periodic.
It looks very natural to conjecture that if only $T$
is an arbitrary expanding polynomial in the sense of Complex
Dynamics \cite {EL} then its
balanced measure  produces a
limit periodic Jacobi matrix. Several research groups attacked this
problem (in full generality)  but failed. Even the case
of the quadratic polynomial with $C>2$ is still open.


Recall some properties of Jacobi matrices. Let $J$
   be  a Jacobi matrix, $J^*=J$, acting in $\CC^d$
   or $l^2(\ZZ_+)$.
Under the assumption $p_k\not=0$ the vector $|0\rangle$ of the
standard basis
   is cyclic for $J$.
   The resolvent function is a function of the
form
\begin{equation}\label{3a}
r(z)=\left<0\left|(J-z)^{-1} \right|0\right>.
\end{equation}
It has positive imaginary part in the upper half plane and hence
possesses the representation
\begin{equation}\label{3b}
   r(z)=\int\frac{d\sigma}{\lambda-z}=
\left<\one
\left|(\lambda-z)^{-1} \right|\one\right>_{L_{d\sigma}^2}.
\end{equation}
where $\lambda\cdot$ is the
operator multiplication by the independent variable in
$L^2_{d\sigma}$ and $\one$ is the function that equals one
identically. 
Formulas \eqref{3a} and \eqref{3b} give one to one correspondence
between triples 
$\{L^2_{d\sigma},\lambda\cdot, \one\}$ and
$\{l^2(\ZZ_+), J, |0\rangle\}$ or
$\{\CC^d, J, |0\rangle\}$, respectively, in the finite
dimensional case.
To restore $J$ starting from the nonnegative measure $\sigma$ one uses
\eqref{r1}.

Our first object is the following
\begin{conjecture}   \label{p}
Let $T(z)$ be an expanding polynomial of degree $d$ with a real
Julia set $E$, $ E\subset [-\xi,\xi]$,
$T^{-1}:[-\xi,\xi]\to[-\xi,\xi]$. Define $J=J(x)$ by
\begin{equation}\label{p0}
\left<0\left|(z-J(x))^{-1}
\right|0\right>=\frac{T'(z)/d}{T(z)-x},\quad x \in[-\xi,\xi].
\end{equation}
   Respectively $J_n(x)$ is associated with an iteration 
   $T_n=T^{\circ n}$, $\deg T_n=d_n$.
  Then for every $\epsilon$ there exists
$n$ such that
\begin{equation}\label{p1}
||J_n(x)-J_n(0)||\le \epsilon.
\end{equation}
\end{conjecture}
\noindent
Note that eigenvalues of $J_n(x)$ and $J_n(0)$ are close,
so the non trivial part  deals with eigenvectors.

Let us explain how this conjecture is related
to the general one. If 
   $\mu$ is the balanced measure on
   $E$, then the resolvent of
$J=J(\mu)$ satisfies to the following Renormalization Equation
\begin{equation}\label{p100}
V^*(z-J)^{-1}V=(T(z)- J)^{-1}T'(z)/d,
\end{equation}
where $V|k\rangle=|kd\rangle$. Let us include $J$ into a chain
$\{J_n(t)\}_{ t\in [0, 1]}$ defined by
\begin{equation*}
V_n^*(z-J_n(t))^{-1}V_n=(T_n(z)-t J)^{-1}T_n'(z)/d_n
\end{equation*}
(compare the last equation with \eqref{p0}).
Then the main goal is to show that 
$$
||J_n(1)-J_n(0)||\le \epsilon
\quad
\text {for}\   n>n_0,
$$ 
since it would imply immediately   that
   $J(\mu)$ is limit periodic. Thus to prove
   Conjecture \ref{p} is a good
model problem on the way to prove limit periodicity of
$J(\mu)$.

The following approach looks very natural: to get \eqref{p1}
we have to estimate $J'(x)$. The given
derivative has a special representation
\begin{equation}\label{p101}
\frac{dJ(x)}{dx}=F(J)+[G,J]
\end{equation}
with $F(J)=\{T'(J)\}^{-1}$.
It is a certain flow on Jacobi matrices that
in a sense is dual to the well--known Toda flow.
We call it FG flow \cite{FG}
(see Sect. 2). 
The first term at the right hand side in
\eqref{p101} is small due to the characteristic
property of expanding polynomials:
$|T'_n(x)|\ge A c^n$,
$x\in E$, with $A>0$, $c>1$. It appears that the estimation
we get for $G$ is not enough to state that
the commutator $[G,J]$ is sufficiently small
(Proposition \ref{2.5}). However on this way we found
 quite designing formulas
 and connections with so many objects
 (the Hilbert transform, the Schwarz derivative,
 the Ruelle and Laplace operators)
 that, we sure, they are of independent
interest. 

In the framework of this approach, initiated in \cite{SYU},
we managed to prove the following theorem that
partially confirms the main hypothesis.

\begin{theorem}\label{3.5}
    Let $J$ be the Jacobi matrix
associated with iterations of an expanding polynomial $T$. Then
for every $\epsilon$
    there exists $n$ such that
\begin{equation}\label{14}
|p_{k+sd_n^2}-p_k|\le\epsilon,\quad |q_{k+sd_n^2}-q_k|\le\epsilon,
\end{equation}
for all $s\ge 0$ and $k=1,2,...,d_n$.

\end{theorem}
\noindent
Note that actually our goal is to prove \eqref{14} when
$k=1,2,...,d_n^2$.
A proof of the theorem is given in Sect. 3.

\medskip
\noindent
{\bf Acknowledgment.} We wish to thank Misha Shapiro who called our
attention to the results of \cite{FG}.

\section{FG flow}

\subsection{Definition} Let $J:\CC^d\to \CC^d$.
Consider the resolvent function
\begin{equation}
\left<0\left|(z-J)^{-1} \right|0\right>=
\sum_{k=1}^d\frac{\sigma_k}{z-\lambda_k}.
\end{equation}
    Under the
Toda flow the spectrum is stable $\lambda_k=\text{Const}$ but
masses vary with time $\sigma_k=\sigma_k(t)$. In FG flow case
$\lambda_k=\lambda_k(t)$ but $\sigma_k=\text{Const}$. Moreover, in
our case \eqref{p0} time is $x$, $\sigma_k=1/d$ and $T(\lambda_k(x))=x$.
Recall $T(z)$ is an expanding polynomial of degree $d$ with a real
Julia set $E$, $ E\subset [-\xi,\xi]$,
$T^{-1}:[-\xi,\xi]\to[-\xi,\xi]$.

We want to get a differential equation on $J$. Let $\fB$ be a
unitary matrix such that
$$
J\fB=\fB\Lambda,
$$
where $\Lambda=\text{diag}\{\lambda_k\}$. Since we can choose
$$
\lambda_1(x)< \lambda_2(x)<...< \lambda_d(x)
$$
   that holds for all $x$, $\fB$ essentially is well defined.
We put
\begin{equation*}
\fB= \frac 1{\sqrt{d}}
\begin{bmatrix}
P_0(\lambda_1)&...&P_0(\lambda_d)\\
\vdots& &\vdots\\
P_{d-1}(\lambda_1)&...&P_{d-1}(\lambda_d)
\end{bmatrix},
\end{equation*}
where $P_k(z)$ is the orthonormal polynomial.

    We differentiate $J$ with respect to $x$
$$
\dot J=\fB\dot\Lambda\fB^{-1}+ \dot\fB\Lambda\fB^{-1}-
\fB\Lambda\fB^{-1} \dot\fB \fB^{-1} =F+GJ-JG,
$$
where $F:=\fB\dot\Lambda\fB^{-1}$, $G:=\dot\fB\fB^{-1}$. By the
definition $F=f(J)$ with $f(\lambda_k)=\dot\lambda_k$. Thus
$F=T'(J)^{-1}$. The next step  is to determine $G$.

Note some evident facts. $G$ is skew--symmetric and $\langle 0| G=
0$, so $G|0\rangle=0$. Also it is easy to show, say by induction,
that
$$
\frac{d}{dx}J^n= n J^{n-1} F+ G J^n-J^n G.
$$
Finally, since $P_k(J)|0\rangle=|k\rangle$ and
$$
\frac{d}{dx}P_k(J)-\frac{\partial} {\partial x}P_k(J)= F P'_k(J)+G
P_k(J)-P_k(J) G
$$
we get
\begin{equation}\label{16}
-\frac{\partial} {\partial x}P_k(J)|0\rangle= F P'_k(J)|0\rangle+G
P_k(J)|0\rangle.
\end{equation}

Let $G_+$ be a lower triangle matrix with zeros on the main
diagonal such that $G=G_+-G_+^*$. Then \eqref{16} implies
\begin{equation}\label{17}
G_+|k\rangle=G_+P_k(J)|0\rangle=
   -(FP'_k(J)|0\rangle)_{+}^{(k)}.
\end{equation}
Here
    $h_{+}^{(k)}$ means that in a vector
$h=\{h_j\}_{j=0}^{d-1}$
    we have to replace
     all coordinates $h_j$, $0\le j\le k$, by zeros.

Let us rewrite \eqref{17} in other words.
   Define an operator $D$ by
$$
D|k\rangle=DP_k(J)|0\rangle:= P'_k(J)|0\rangle.
$$
Then
$$
G_+=-(FD)_+.
$$

\smallskip

It is easy to check using the functional representation in
$L^2_{d\sigma}$ that
\begin{equation}\label{18}
D J- J D= I- |(p_dP_d)'\rangle\langle P_{d-1}|,
\end{equation}
where $p_d P_d(\lambda)$ is defined by \eqref{r1}. Note that $p_d
P_d(z)=0$ in $L^2_{d\sigma}$, that is it has the
same roots
$\{\lambda_k(x)\}$ as ${T(z)-x}$. Thus
   $p_d P_d(z)=C(T(z)-x)$ and  $(p_d P_d)'(z)=CT'(z)$.

\begin{definition} FG flow is given by a differential equation
of the form
\begin{equation}\label{19}
\dot J=F+G J-J G
\end{equation}
with $F=f(J)$ and $G=G_+ - G_+^*$, where $G_+=-(FD)_+$ and $D$ is
an (upper triangle) matrix such that commutant $[D,J]$ equals the
unity matrix up to a one dimensional perturbation \cite{FG}.
\end{definition}

\subsection{$(FD)$ as a Hilbert transform}
\begin{lemma} The matrix of the operator $(FD)$ with respect to
the basis of eigenvectors of $J$ has the form
\begin{equation}
\begin{bmatrix}
\frac 1 2\frac{T''(\lambda_1)}{T'(\lambda_1)}& \hdots&
\frac 1{\lambda_1-\lambda_d}\\
\vdots & & \vdots\\
\frac 1{\lambda_d-\lambda_1}&\hdots &\frac 1
2\frac{T''(\lambda_d)}{T'(\lambda_d)}
\end{bmatrix}
\begin{bmatrix}
\frac{1}{T'(\lambda_1)}& & \\
&\ddots & \\
& & \frac{1}{T'(\lambda_d)}
\end{bmatrix}.
\end{equation}
\end{lemma}

\begin{proof}
Let us evaluate $D$ in the basis of eigenvectors of $J$. In this
basis
$$
|P(x)\rangle\to\frac 1{\sqrt{d}}\begin{bmatrix} P(\lambda_1)\\
\vdots\\P(\lambda_d)\end{bmatrix}.
$$
   As we know $(p_d P_d)'(\lambda_k)=CT'(\lambda_k)$.
   Taking into account that now $J$ is diagonal we conclude that 
   the diagonal
   entries  of $DJ-JD$ are zeros. Therefore
$P_{d-1}(\lambda_j)= {d}/\{C T'(\lambda_j)\}$. Thus the right hand
side of \eqref{18} is of the form
$$
I-\begin{bmatrix}T'(\lambda_1)\\
\vdots\\T'(\lambda_d)\end{bmatrix}
\begin{bmatrix}\frac 1{T'(\lambda_1)},& \hdots,&
\frac 1{T'(\lambda_d)}\end{bmatrix}.
$$
Referring again to a diagonal form of $J$ we solve \eqref{18} and
get
\begin{equation}\label{a13}
D_{ij}= \frac 1{\lambda_i-\lambda_j}
\frac{T'(\lambda_i)}{T'(\lambda_j)}, \quad i\not=j.
\end{equation}

To find diagonal entries $D_{ii}$ we have to use $D|0\rangle=0$.
Since
\begin{equation*}
|0\rangle\to\one\to\frac 1{\sqrt{d}}\begin{bmatrix} 1\\
\vdots\\ 1\end{bmatrix}
\end{equation*}
we get
$$
D_{ii}=T'(\lambda_i)\sum_{k\not=i}\frac 1{T'(\lambda_k)}\frac
1{\lambda_k-\lambda_i}.
$$
Note that
$$
\sum_{k=1}^d\frac 1{T'(\lambda_k)}\frac 1{\lambda_k-z}=-\frac
1{T(z)-x}.
$$
Therefore
\begin{equation}\label{a14}
\begin{split}
   \frac{D_{ii}}{T'(\lambda_i)}=& \lim_{z\to\lambda_i}
\left\{-\frac 1{T'(\lambda_i)}\frac 1{\lambda_i-z}-\frac
1{T(z)-x}\right\} \\
=& \lim_{z\to\lambda_i}\frac{ \frac{T(z)-x}{z-\lambda_i}\frac
1{T'(\lambda_i)}-1}{T(z)-x}= \frac{\frac 1 2
T''(\lambda_i)}{(T'(\lambda_i))^2}.
\end{split}
\end{equation}
Thus \eqref{a13} and \eqref{a14} finish the proof.
\end{proof}

\subsection{Trace of $(FD)^*(FD)$}

\begin{lemma}Let $L_2$ be a Ruelle operator of the form 
\begin{equation}\label{l2}
L_2g(x)=\frac 1
d\sum_{Ty=x}\left(\frac{g}{{T'}^2}\right)(y)
\end{equation}
and let $S(T)$ be 
the Schwarz derivative of $T$, $S(T)=\frac{T'''}{T'}-\frac 3
2\left(\frac{T''}{T'}\right)^2$ . Then
$$
\frac 1 d\tr\{(FD)^*(FD)\}= -\frac 1 3 L_2 \{S(T)\}.
$$

\end{lemma}

\begin{proof}
First we simplify
\begin{equation}\label{a15}
   u_i=\sum_{k\not=i}\frac{1}{(\lambda_i-\lambda_k)^2}=
\lim_{z\to\lambda_i} \left\{
\sum_{T(\lambda)=x}\frac{1}{(z-\lambda)^2}-
\frac{1}{(z-\lambda_i)^2}\right\}.
\end{equation}
Note that
$$
\sum_{T(\lambda)=x}\frac{1}{(z-\lambda)}=\frac{T'(z)}{T(z)-x}.
$$
That is
$$
\sum_{T(\lambda)=x}\frac{1}{(z-\lambda)^2}=
\frac{{T'}^2(z)-T''(z)(T(z)-x)}{(T(z)-x)^2}.
$$
So passing in a usual way to the limit in \eqref{a15} we get
\begin{equation}
u_i=\frac{\frac 1 2(T'')^2-\frac 2 3 T' T'''}
{2(T')^2}(\lambda_i).
\end{equation}
This means that a diagonal entry of the operator $(FD)^*(FD)$ with
respect to the basis of eigenvectors of $J$ has the form
$$
\frac{1}{{T'}^2}\left\{\left(\frac{1}2\frac{T''}{T'}\right)^2(\lambda_i)+u_i\right\}=
\frac{1}{{T'}^2} \frac{(T'')^2-\frac 2 3 T' T'''}
{2(T')^2}(\lambda_i)=-\frac 1
3\left(\frac{S(T)}{{T'}^2}\right)(\lambda_i).
$$

\end{proof}

   Naturally, in the same way we
can find off diagonal entries of the matrix of
the operator $(FD)^*(FD)$.
\begin{lemma}\label{2.4} For $i\not=j$
\begin{equation*}
\{(FD)^*(FD)\}_{ij}=\frac 1{T'(\lambda_i)}
\frac 2{(\lambda_i-\lambda_j)^2}
\frac 1{T'(\lambda_i)}.
\end{equation*}
\end{lemma}

We would consider 
$\Delta:=(FDF^{-1})^*(FDF^{-1})$ as a counterpart of Laplacian due to the
following proposition.
\begin{corollary}
$\Delta$ is a
positive operator that satisfies
$$[J,[J,\Delta]]=2d|0\rangle\langle 0|-2.$$
\end{corollary}

\begin{proof} See Lemma \ref{2.4}.
\end{proof}

Our plan to estimate $[G,J]$ in \eqref{p101}
was based on the conjecture
 $||(FD)_n||\sim
\kappa^n$ with $\kappa<1$ (typically everything that goes to zero
in the subject goes to zero as a geometric progression). 
Since $(G_n)_+=-(FD)_{n +}$ 
that would give
an estimation on $G$:
$$
||(G_n)_+||\sim\kappa^n n\log d,
$$
and we are done. However
the following proposition  shows that
$||(FD)_n||\not\to 0$.

\begin{proposition}\label{2.5} There exists the limit
\begin{equation}\label{a17}
\lim_{n\to\infty}\frac 1 {d^n}\tr\{(FD)_n^*(FD)_n\}=-\frac 1 3
(I-L_2)^{-1}L_2 S(T).
\end{equation}
\end{proposition}

\begin{proof}
Let us use the Chain Rule for the Schwarz derivative
$$
S(T_{n+1})=S(T_n)\circ T {T'}^2+S(T).
$$
    Since  $L_2\{g\circ T {T'}^2\}=g$ holds for
every function $g$, we have
   $L_2^{n+1}\{S(T_n)\circ T
{T'}^2\}=L_2^n S(T_n)$ and therefore
\begin{equation}\label{a18}
L^{n+1}_2 S(T_{n+1})=L^{n}_2 S(T_n)+L^{n+1}_2 S(T)= L_2
S(T)+...+L^{n+1}_2 S(T).
\end{equation}
The spectral radius of $L_2$ less than $1/d^2$  (see Lemma
\ref{2.1}). So
\eqref{a18} completes the proof.
\end{proof}

\begin{proof}[Remark]
We still believe in the limit periodic
property of
$J(\mu)$. Recall that we have to estimate not $(FD)_n$ itself but the
commutator $[G_n,J_n]$.
Probably it worth to mention that the right
hand side of the commutant identity for
$(FD)_n$,
$$
(FD)_n J_n-J_n(FD)_n= F_n-d^n |0\rangle\langle 0| F_n
$$
goes to zero in norm (it's again Lemma  \ref{2.1}). That is
asymptotically $(FD)_n$ and
$J_n$ commute.

\end{proof}

\section{Partial result in the right direction}

\subsection{
Renormalization equation} Let
$$
Lg(x)=\frac 1 d\sum_{Ty=x}g(y)
$$
be a Ruelle
operator associated with
 an expanding polynomial $T(z)$. If $\tilde J$ is the Jacobi
matrix associated with a measure
$\tilde\sigma$ supported on $E$, $\tilde J:= \tilde J(\tilde\sigma)$,
then the Renormalization Equation
\begin{equation}\label{t01}
V^*(z-J)^{-1}V=(T(z)-\tilde J)^{-1}T'(z)/d,
\quad V|k\rangle=|kd\rangle,
\end{equation}
has a unique solution
$ J:= J(\sigma)$,  where $\sigma:=
L^*(\tilde\sigma)$ \cite{BGH}, \cite{LSYU}.
 It follows basically from the identity
$$
\left( L\frac{1}{z-y}(g\circ T)(y)\right)(x)=
\frac{T'(z)/d}{T(z)-x}g(x)
$$
and the functional representations of both operators
in $L^2_{d\sigma}$ and $L^2_{d\tilde\sigma}$
respectively. Note that \eqref{t01} becomes \eqref{p100}
if $\tilde\sigma=\mu$, since for the balanced measure we have
$\mu=
L^*(\mu)$.

\begin{lemma}\label{1.1} Let $J^{(s)}$ be the
$s$-th $d\times d$ block of the matrix $J$, that is
\begin{equation}
J^{(s)}=
\begin{bmatrix}
q_{sd}&p_{sd+1}& & & \\
p_{sd+1}&q_{sd+1}&p_{sd+2}& & \\
    &\ddots&\ddots&\ddots & \\
& & p_{sd+d-2}&q_{sd+d-2}&p_{sd+d-1} \\
&  & &p_{sd+d-1}&q_{sd+d-1}
\end{bmatrix}.
\end{equation}
Then its resolvent function is of the form
\begin{equation}\label{3}
\left<0\left|(z-J^{(s)})^{-1} \right|0\right>=
\frac{T'(z)/d}{T^{(s)}(z)}.
\end{equation}
Moreover at the critical points $\{c: T'(c)=0\}$ the following
decomposition in a continued fraction holds true
\begin{equation}
T^{(s)}(c)=T(c)-\tilde q_s- \frac{\tilde p^2_s} {T(c)-\tilde
q_{s-1}-...}.
\end{equation}
\end{lemma}

\begin{proof}
We write $J$ as a $d\times d$ block matrix (each block is of
infinite size):
\begin{equation}\label{5}
J=
\begin{bmatrix}
\cQ_{0}&\cP_{1}& & &S_+\cP_d  \\
\cP_{1}&\cQ_{1}&\cP_{2}& & \\
    &\ddots&\ddots&\ddots & \\
& & \cP_{d-2}&\cQ_{d-2}&\cP_{d-1} \\
\cP_d S^*_+&  & &\cP_{d-1}&\cQ_{d-1}
\end{bmatrix}.
\end{equation}
Here $\cP_k$ (respectively $\cQ_k$) is a diagonal matrix
$\cP_k=\text{diaq}\{p_{k+sd}\}_{s\ge 0}$ and $S_+$ is the
one--sided shift.  In this case $V^*$ is the projection on the
first block--component.

Using this representation and being well
    known identity for  block matrices
\begin{equation*}
\begin{bmatrix}
A&B\\
C&D
\end{bmatrix}^{-1}=
\begin{bmatrix}
(A-B D^{-1}C)^{-1}&*\\
*&*
\end{bmatrix},
\end{equation*}
we get
\begin{equation}\label{6}
\frac{T(z)-\tilde J}{T'(z)/d}= z-\cQ_0-
\begin{bmatrix}\cP_1,&...,&S_+\cP_d
\end{bmatrix}
\{z-J_{1}\}^{-1}
\begin{bmatrix}
\cP_1\\
\vdots\\
\cP_d S_+^*
\end{bmatrix},
\end{equation}
where $J_{1}$ is the matrix that we obtain from $J$ by deleting
the first block--row and the first block--column in \eqref{5}.
Note that in $(z-J_{1})$ each block is a diagonal matrix that's
why we can easily get an inverse matrix in terms of orthogonal
polynomials.

Let us introduce the following notations: everything related to
$J^{(s)}$ has superscript $s$. For instance: $p_k^{(s)}=p_{sd+k}$,
$1\le k\le d$, respectively $P^{(s)}_d$ and $Q^{(s)}_d$ mean
orthonormal polynomials of the first and second kind. In this
terms equation \eqref{6} is equivalent to the two series of scalar
relations corresponding to the diagonal and off diagonal entries
\begin{equation}\label{7}
\frac{T(z)-\tilde q_{s+1}}{T'(z)/d}= \frac{P_d^{(s+1)}(z)}
{Q^{(s+1)}_d(z)}- p_{ds}^2\frac{Q_{d-1}^{(s)}(z)/p_{ds}}
{Q^{(s)}_d(z)}
\end{equation}
     and
\begin{equation}\label{8}
\frac{\tilde p_{s+1}}{T'(z)/d}= \frac{p_1^{(s)} ... p_d^{(s)}}
{z^{d-1}+...}=\frac 1{Q^{(s)}_d(z)}.
\end{equation}
We have to remind (see \eqref {3} and \eqref{8}) that
\begin{equation*}
\frac{Q_d^{(s)}(z)} {P^{(s)}_d(z)} =\frac{z^{d-1}+...}{z^d+...} =
\frac{T'(z)/d} {T^{(s)}(z)}.
\end{equation*}

    Now, due to the
Wronskian identity, if $T'(c)=0$ then
\begin{equation}\label{9}
-p_{ds}{Q^{(s)}_{d-1}(c)} =\frac 1{P^{(s)}_d(c)}.
\end{equation}
So, combining \eqref{7}, \eqref{8} and \eqref{9} we get the
recurrence relation
\begin{equation}
T(c)-\tilde q_{s+1}= T^{(s+1)}(c)+\frac{\tilde p^2_{s+1}}
{T^{(s)}(c)}
\end{equation}
with initial data
$$
T^{(0)}(c)=T(c)-\tilde q_0.
$$

\end{proof}

\subsection{ $p_{s d_n}$ are exponentially small}

\begin{lemma}\label{2.1}
    Let $J$ be the Jacobi matrix
associated with iterations $\{T_n\}_{n\ge 1}$ of an expanding
polynomial $T$, that is $J=J(\mu)$ where $L^*\mu=\mu$. Then
\begin{equation}\label{11}
C_-(\rho d)^n p_s\le p_{s d_n}\le C_+(\rho d)^n p_s
\end{equation}
with $C_\pm>0$ and $0<\rho<1/d$.
\end{lemma}

\begin{proof}
We recall that $P_{sd}=P_s\circ T$ and $Q_{sd}= (T'/d) Q_s\circ T
$ \cite{BGH}, \cite{LSYU}. We use an interpolation formula
\begin{equation}\label{12}
\int R\,d\mu= \sum_{y: P_{sd}(y)=0} R(y)\frac{Q_{sd}}
{P'_{sd}}(y),\quad \deg R<sd,
\end{equation}
and the Wronskian identity
\begin{equation}\label{13}
p_{sd}\{P_{sd-1} Q_{sd}- Q_{sd-1} P_{sd}\}=1.
\end{equation}
Substituting \eqref{13} in \eqref{12} we obtain
\begin{equation*}
p_{sd}^2=\int \{p_{sd} P_{sd-1}\}^2\,d\mu= \sum_{y:
P_{sd}(y)=0}\{p_{sd} P_{sd-1}(y)\}^2 \frac{Q_{sd}}{P'_{sd}}(y)=
\sum_{y: P_{sd}(y)=0} \frac{1}{(Q_{sd} P'_{sd})(y)}.
\end{equation*}
Therefore,
\begin{equation*}
\begin{split}
p_{sd}^2=& \sum_{x: P_{s}(x)=0} \sum_{y: T(y)=x}
\frac{1}{(({T'}^2/d)(Q_{s} P'_{s})\circ T )(y)}\\=& \sum_{x:
P_{s}(x)=0} \frac{1}{(Q_{s} P'_{s})(x)}\left\{ \frac 1 d\sum_{y:
T(y)=x}\frac {d^2} {{T'}^2(y)}\right\}\\=& \sum_{x: P_{s}(x)=0}
\{p_{s} P_{s-1}(x)\}^2 \frac{Q_{s}}{P'_{s}}(x)\left\{ \frac 1
d\sum_{y: T(y)=x}\frac {d^2} {{T'}^2(y)}\right\}.
\end{split}
\end{equation*}
Now we use the Ruelle version of the Perron--Frobenius theorem
\cite{EL}, \cite{B}
with respect to
$L_2$ \eqref{l2}. 
According to this theorem
$$
\frac 1{\rho^{2n}} L^n_2 g\to h(x)\int g\,d\nu,
$$
uniformly on $x$ with a certain continuous function $h>0$ and
positive measure $\nu$; $\rho^2$ is the spectral radius of $L_2$.
Combining this with the interpolation formula we get both--sided
estimate \eqref{11}.

We only have to show that $(\rho d)<1$.
Let $b(z)$ be the complex Green's function of the domain 
$\overline \CC \setminus E$ with respect to infinity. 
Consider the sequence of functions
$\{f_n\}_{n\ge 1}$, where $f_n(z):=(b^{d_n-1} P_{d_n-1})(z)$.
It is a multiple--valued function in the domain
$\overline \CC \setminus E$ with a single--valued modulus
which has a harmonic majorant $u_n(z)$: $|f_n(z)|^2\le u_n(z)$.
Moreover, $u_n(\infty)=||P_{d_n-1}||^2_{L^2_{d\mu}}=1$.
We claim that $f_n$ should go to zero pointwise.
If not then we can find a subsequence $\{f_{n_k}\}$ that
converges to a non trivial function $f$. However, in this case,
$(bf)(z)$ is a non trivial single valued in $\overline \CC \setminus E$
function,
$|(bf)(z)|^2$ has a harmonic majorant and
$(bf)(\infty)=0$. This contradicts to the well-known fact
that {\it analytic} capacity (that is the Lebesgue measure in this case)
of $E$ is zero.

 Therefore the sequence converges to zero.
In particular
$$
(b^{d_n-1} P_{d_n-1})(\infty)=\frac{1}{p_1...p_{d_n-1}}=
\frac{p_{d_n}}{p_1}\to 0,\quad  n\to \infty.
$$
But $\frac{p_{d_n}}{p_1}\sim (\rho d)^n$, thus $(\rho d)
<1$.
\end{proof}

\begin{remark}
Let us mention here that $q_{sd}=q_0$ since
$$
q_{sd}=\int y P^2_{sd}\,d\mu =\int y P_s^2\circ T\,dL^*\mu =\int
(Ly) P_s^2\,d\mu
$$
and $Ly=q_0$.
\end{remark}

\subsection{The result}
First we prove  (undoubtedly well--known and  simple)
\begin{lemma}\label{3.1}
Assume that two measures $\sigma$ and $\tilde\sigma$
are mutually absolutely continuous. Moreover,
 $d\tilde\sigma
=f\,d\sigma$ and 
$1-\epsilon\le f\le (1-\epsilon)^{-1}$.
    Let us associate with these measures Jacobi matrices
$J= J(\sigma)$,
$\tilde J= J(\tilde\sigma)$.
 Then for their coefficients we have
$$
|\tilde p_s-p_s|\le \frac\epsilon{1-\epsilon}|| J||.
$$
\end{lemma}

\begin{proof}
Let us use an extreme property of orthogonal polynomials,
\begin{equation*}
\begin{split}
\tilde p_1^2...\tilde p_s^2=& \int\tilde p_1^2...\tilde p_s^2
\tilde P_s^2\,d\tilde\sigma \ge
(1-\epsilon)\int\{z^s+...\}^2\,d \sigma\\
\ge&(1-\epsilon) \inf_{\{P=z^s+...\}}\int P^2 \,d \sigma
=(1-\epsilon) p_1^2... p_s^2.
\end{split}
\end{equation*}
Similarly
\begin{equation*}
p_1^2... p_{s-1}^2
    \ge
(1-\epsilon)\tilde p_1^2...\tilde p_{s-1}^2.
\end{equation*}
Therefore
$$
\frac 1{(1-\epsilon)^2}p_s^2\ge \tilde p_s^2\ge(1-\epsilon)^2p_s^2
$$
    and hence
$$
-\epsilon p_s\le \tilde p_s-p_s\le \frac\epsilon{1-\epsilon} p_s.
$$

\end{proof}

Now, we are in position to prove Theorem \ref{3.5}.

\begin{proof}
As it follows from Lemma \ref{1.1}
\begin{equation*}
T^{(s)}(c)=T(c)-q_s-p_s^2\int \frac{d\nu^{(s)}(x)} {T(c)-x}.
\end{equation*}
Here $\nu^{(s)}$ is a discrete measure such that
$\text{supp}\{\nu^{(s)}\}\subset[-\xi,\xi]$,
$\nu^{(s)}([-\xi,\xi])=1$, recall that $[-\xi,\xi]$ is the smallest
interval containing the Julia set $E$. In particular,
\begin{equation*}
T^{(sd_n)}(c)=T(c)-q_0-p_{sd_n}^2\int \frac{d\nu^{(sd_n)}(x)}
{T(c)-x}.
\end{equation*}
Now, since
$$
\text{dist}_{\{c:T'(c)=0\}}\{T(c),[-\xi,\xi]\} =\delta>0,
$$
for every $\epsilon>0$ there exists $n$ such that
$$
(1-\epsilon)\le \frac{T^{(sd_n)}(c)}{T(c)-q_0}
\le(1-\epsilon)^{-1}
$$
    (here we  used Lemma \ref{2.1}).
Recall \eqref{3}, so, Lemma \ref{3.1} completes
the proof.

\end{proof}

\bibliographystyle{amsplain}

\begin{thebibliography}{1}

\bibitem{AS} J. Avron and B. Simon,
\textit{Singular continuous spectrum for a class
 of almost periodic Jacobi matrices}
 Bull. AMS {\bf 6} (1982), 81--85.

\bibitem{BGH}
 M. F. Barnsley, J. S. Geronimo,  
A. N. Harrington,
\textit{Almost periodic Jacobi matrices associated 
with Julia sets for polynomials},
Comm. Math. Phys. {\bf 99} (1985), no. 3, 303--317.

\bibitem{BBM} J.  Bellissard,  D. Bessis, P. Moussa,  
\textit {Chaotic
states of almost periodic Schr\"odinger operators}, Phys. Rev. Lett. 
{\bf 49}
(1982), no. 10, 701--704.


\bibitem{BS}
J. Bellissard, B. Simon,  
\textit {Cantor spectrum for the
almost Mathieu equation}
J. Funct. Anal. {\bf 48} (1982), no. 3, 408--419.

\bibitem{B}
R. Bowen,
\textit{
Equilibrium states and the ergodic theory
of Anosov diffeomorphisms}
Berlin, Heidelberg, New York: SPringer 1975.

\bibitem{FG}
 L. Faybusovich and M. Gekhtman,
 \textit{ Poisson
brackets on rational functions and multi-Hamiltonian structure for
integrable lattices}, Phys. Lett. A {\bf 272} (2000), no. 4, 236--244.

\bibitem{EL}
A. Eremenko and M. Lyubich,
\textit {The dynamics of analytic
transformations}. (Russian) Algebra i Analiz 1 (1989), no. 3, 1--70;
translation in Leningrad Math. J. {\bf 1} (1990), no. 3, 563--634.

\bibitem{LSYU}
G. Levin, M. Sodin, P. Yuditskii, 
\textit{
A Ruelle operator for a real
Julia set}, Comm. Math. Phys. {\bf 141} (1991), no. 1, 119--132. Ê

\bibitem{SYU}
M. Sodin, P. Yuditskii, 
\textit{
The limit-periodic finite-difference operator
on
$l\sp 2(\bold Z)$ associated with iterations of quadratic polynomials},
J. Statist. Phys. {\bf 60} (1990), no. 5-6, 863--873.

\end{thebibliography}

\end{document}